\documentclass[14pt,a4paper, preprint]{article}

\usepackage[dvipsnames]{xcolor}
\usepackage{color}
\usepackage[utf8]{inputenc}
\usepackage{amsmath}
\usepackage{amssymb}
\usepackage{float}
\usepackage{graphicx}
\usepackage{cite}
\usepackage{enumerate}
\usepackage[colorinlistoftodos]{todonotes}
\usepackage{appendix}
\usepackage{upgreek}
\usepackage{color, colortbl}
\usepackage{etoolbox}
\usepackage{hyperref}
\usepackage{dcolumn}
\usepackage{fullpage}
\usepackage{relsize}
\usepackage{url}
\usepackage{textcomp}
\usepackage{lmodern}
\usepackage[T1]{fontenc}
\usepackage{authblk}
\usepackage{xfrac}
\usepackage{graphicx}
\usepackage{mathrsfs}
\usepackage{dcolumn}
\usepackage{float}
\usepackage{geometry}
\usepackage{multicol}
\usepackage{tikz}
\usepackage{pgfplots}
\usepackage{etoolbox}
\usepackage{enumerate}
\usepackage{cite}
\usepackage{subfiles}
\usepackage{lettrine}
\usepackage{import}
\usepackage{pdfpages}
\usepackage{bm}

\usepackage{algorithm}
\usepackage{algpseudocode}

\usepackage{tikz}
\usetikzlibrary{fit}

\hypersetup{
    colorlinks=true,
    linkcolor=blue,
    filecolor=magenta,      
    urlcolor=blue,
    pdfpagemode=FullScreen,
    }

\definecolor{Gray}{gray}{0.9}
\definecolor{my_emerald}{RGB}{1,99,52}
\definecolor{darkseagreen}{rgb}{0.56, 0.74, 0.56}

\makeatletter
\renewcommand*\env@matrix[1][*\c@MaxMatrixCols c]{%
  \hskip -\arraycolsep
  \let\@ifnextchar\new@ifnextchar
  \array{#1}}
\makeatother

\def\env@matrix{\hskip -\arraycolsep
  \let\@ifnextchar\new@ifnextchar
  \array{*\c@MaxMatrixCols c}}

\newenvironment{amatrix}[1]{%
  \left[\begin{array}{@{}*{#1}{c}|c@{}}
}{%
  \end{array}\right]
}

\usepackage[inline]{trackchanges}

\addeditor{MR}
\addeditor{RA}
\addeditor{AL}

\title{Designing Problems for Improved Instruction and Learning - Linear Algebra}

\author[1]{R. H. Allaire}
    \affil[1]{Department of Mathematical Sciences, United States Military Academy, West Point, NY, 10996}
\author[1]{M. Reynolds}

\author[1]{A. Lee}

\date{}

\begin{document}
\maketitle

\begin{abstract}
    One of the grand challenges of Mathematics instruction is to provide students with problems that are both accessible and have a reasonably elegant solution. Instructors commonly resort to resources like course textbooks, online-learning platforms, or other automated problem-generating software to select problems for exams and assignments. However, reliance on such tools may result in limited control over problem parameters, potentially yielding intricate solutions that impede students' understanding. This article centers on Linear Algebra, wherein we devise algorithms for reverse engineering matrices of integers with integer outcomes through operations such as the inverse, LU decomposition, and QR decomposition. The focus is on empowering instructors to manipulate matrix properties deliberately, ensuring the creation of problems that enrich instruction and foster student confidence. The intellectual endeavor of reverse engineering such problems, grounded in both theory and matrix properties, proves mutually beneficial for both students and instructors alike.
\end{abstract}

\section{Introduction}
Suppose you were a math instructor guiding students through the concept of solving a system of linear equations and you pose the following problem, either in class, on an online homework assignment, or on an exam. 
\begin{align*}
2x+5y=\underline{16}, \\ 
-x+8y=11. \nonumber
\end{align*}

The solution turns out to be an ``ugly'' $(\frac{73}{21}, \frac{38}{21})$ or $(3.476, 1.809)$, and the interim substitution or row reduction steps contain several fractions and algebraic steps where students can easily make mistakes. If the average student encounters several of these algebraic stumbling blocks along the way, this distracts from the more important solution approach and critical lesson objective of solving a system of linear equations.

Now consider that if you just change just one number in the first problem (the underlined \underline{16} to a 20), the solution would be a ``nice'' solution of $x=5,\, y=2$, with very little algebraic distraction in the way of demonstrating the solution approach. Here we define ``nice'' as an integer solution. 

Likewise, perhaps you inherited an online learning platform to teach how to compute eigenvalues to students and the pre-programmed parameter randomization presents a student with the following problem. 
\begin{align*}
\text{Find all eigenvalues of the matrix } A = \begin{bmatrix} 
    1 & \underline{2} \\
    3 & 5 
    \end{bmatrix}.
\end{align*}
The solution to this eigenvalue problem is unfortunately an ``ugly'' $3 \pm \sqrt{10}$. If the correct answer is not inputted properly into the online learning platform (either using the specified math notation or a decimal approximation with enough significant figures), many students will experience frustrations based on the automated grading parameters of the learning platform. For this eigenvalues problem, if you just change the underlined \underline{2} to a -1, the eigenvalues are ``nice'' integers, $(4, 2)$. 

For students delving into a new math concept like these for the first time, posing a problem with ugly solutions accompanied with ugly interim steps can shake the confidence of students, or worse, push them away from truly comprehending the lesson content. Several studies have shown that learners with higher confidence are more willing to learn, challenge themselves, and have better academic resilience \cite{oro24261, LEE2017269, MARTIN200853}. As these papers indicate, what is most important is building the confidence or self-efficacy of students by teaching first principles and fundamentals first, which contributes to generating a passion for lifelong learning. Building confidence and resilience is paramount to reduce well-studied anxiety with math, which can disrupt cognitive processing by compromising ongoing activity in working memory \cite{doi:10.1111/1467-8721.00196, Baloglu, Hoffman2010}. Using a tree analogy, teachers can branch off to more complex problems later, but they should always aim to start with the simplest root problem and branch out from there. This idea of starting with a simple problem with nice solutions is also relevant to the curating of exam problems, where the foundations of self-efficacy are laid down. 

With respect to developing nice linear algebra problems, there exists some related literature on constructing integer matrices with integer eigenvalues \cite{Towse, https://doi.org/10.48550/arxiv.0712.3060, Martin.2008}. Towse et al. acknowledge the importance of concrete ``toy'' problems as a pedagogical approach to teaching students the topic of diagonalizing matrices. The authors define a matrix A as ``IMIE'' if it is an \underline{I}nteger-Entry \underline{M}atrix with all \underline{I}nteger \underline{E}igenvalues and reference a webtool to construct such matrices. Starting in the mid-1980’s, several papers addressed the construction of integer matrices with special properties \cite{doi:10.1080/00029890.1986.11971837, doi:10.1080/00029890.1984.11971414, 10.2307/2322394}. Others discussed variants to this problem, such as how to generate integer matrices whose inverse only contain integers \cite{10.2307/3026530}. For more advanced topics, others have posed approaches to achieve a fraction-free inverse of a Toeplitz Matrix \cite{Bistritz.2010}.

Similarly, for introductory calculus classes, some have discussed systematic approaches to find cubic polynomials with integer root values \cite{10.2307/2691448}. Buddenhagen et al. also discuss how to form related problems like the box problem (i.e., maximize the volume of an open box, which is constructed from a rectangular piece of metal by cutting four equal squares from the corners and bending up the resulting side tabs). Cuoco et al. calls the problem of constructing a nice problem such as this a ``meta-problem,'' and discusses how to create nice problems that are easier to correct and do not involve messy calculations \cite{doi:10.1080/07468342.2000.11974176}. He describes how classical algebra and number theory can be applied to solving these meta-problems such as the generation of Pythagorean triples and the well-studied diophantine equations (i.e., polynomials with integer coefficients where the solutions are also integer). 

The problems discussed in the literature entail a degree of reverse engineering, where the goal is to provide students with a well-structured matrix or equation. This approach enables them to focus on practicing techniques without being burdened by extraneous and distracting computational tasks.

\section{Methodology} \label{sect:methodology}
In this paper we focus on two central topics in Linear Algebra: (i) matrix decomposition and (ii) solving a linear system, $A\bm{x}=\bm{b}$. More specifically we propose questions such as:
\begin{enumerate}[(i)]
    \item How do we create a matrix $A$ with nice coefficients and controllable properties?
    \item How do we construct a system of equations with integer coefficients that has a unique integer solution, infinitely many solutions, or no solution?
\end{enumerate} 

In Sect.~\ref{sect:matrix_properties} we discuss how to utilize the many factorizations of a matrix to construct a matrix with imposed properties. Then, in Sect.~\ref{sect:linear_systems} we address (ii) and outline algorithms for constructing linear systems with the desired number of solutions. In both Sections \ref{sect:matrix_properties} and \ref{sect:linear_systems} we desire that the coefficients of $A$ are ``nice'' (either in the integer or rational sense), which is an added constraint that makes the construction of such problems difficult.

\subsection{Controlling the Properties of Matrix} \label{sect:matrix_properties}
One may argue elementary Linear Algebra is simply showing numerous factorizations of matrices. For example, Strang outlines in his newest text that Linear Algebra can be viewed as simply trying to factor a matrix into those such as the LU, QR, and CR factorizations \cite{strang2023}. Both diagonalization and singular value decomposition can be similarly viewed as some factorization of a matrix, where the decomposed form holds some level of information. We adopt this approach in this Section, and ask how do we control the properties of a matrix while retaining a ``nice" structure?

\subsubsection{LU Decomposition of Integer Matrices}
Here, we focus on constructing the methods of constructing integer matrices from the $LU$ decomposition and its variants. We ask, how do we utilize $A=LU$, where $L$ is lower-triangular and $U$ is upper-triangular, to create an integer matrix that is (i) full-rank or (ii) rank-deficient? 

Regarding (i), integer-filled $L$ and $U$ lead to a trivially generated matrix $A$ and gives complete user-control to all-elements of $L$ and $U$ (although the convention in many cases is to maintain ones along the diagonal of $L$). A more interesting case is the utilization of non-integer rational numbers in $L$ and $U$. Consider the following $2 \times 2$ example:
\begin{align}
    A= LU=\begin{bmatrix}
        1 & 0 \\
        5/2 & 1
    \end{bmatrix} 
    \begin{bmatrix}
        2 & 4 \\
        0 & 1
    \end{bmatrix}.  \label{eq:rational_2by2_system1}
\end{align}
In Eq.~\eqref{eq:rational_2by2_system1} even though $l_{21}=5/2$, we have that $l_{21}u_{11}$ and $l_{22}u_{12}$ are integers. If either $L$ or $U$ have multiple non-integer rational numbers there is added complexity (although the use of rational numbers along the diagonal of $L$ is non-standard).  
For example, the factorization
\begin{align}
    LU=\begin{bmatrix}
        5/9 & 0 \\
        3/2 & 1
    \end{bmatrix} 
    \begin{bmatrix}
        N \cdot lcm(2,9) & M\cdot lcm(2,9) \\
        0 & 1
    \end{bmatrix}, \label{eq:LU_generalized}
\end{align}
where $N,M\in \mathbb{Z}$, is sufficient to make an integer matrix $A=LU$, since $lcm(2,9)$ divides the first row of $U$. To generalize, we let $\bm{l}^{\rm col}_k=\bm{p}_k/q_k$ represent the $k^{\rm th}$ column of $L$ where $\bm{p}_k \in \mathbb{Q}^n$ and $q_k$ is the least common multiple of the denominators of each of the rational components of $\bm{l}^{\rm col}_k$, and represent the $k^{\rm th}$ row of $U$ by $\bm{u}^{\rm row}_k$. Equation \eqref{eq:LU_generalized} may be then written as:
\begin{align}
    \frac{1}{18}\begin{bmatrix}
        10 \\
        27
    \end{bmatrix} \begin{bmatrix}
        N\cdot lcm(2,9) \\ M \cdot lcm(2,9)
    \end{bmatrix}^T +
    \begin{bmatrix}
        0 \\
        1
    \end{bmatrix} \begin{bmatrix}
        0 \\ 1
    \end{bmatrix}^T = \bm{l}^{\rm col}_1 \otimes \bm{u}^{\rm row}_1 + \bm{l}^{\rm col}_2 \otimes \bm{u}^{\rm row}_2,
\end{align}
where $\otimes$ is the outer product. In the above example $\bm{p}_1= \begin{bmatrix} 10 & 27 \end{bmatrix}^T$, $\bm{p}_2=\begin{bmatrix} 0 & 1 \end{bmatrix}^T$, and it is evident that the least-common multiple of the denominators of $\bm{l}^{\rm col}_1$, $q_1=18$, must divide all of $\bm{u}^{\rm row}_1$. 
For $n \times n$ matrices, the matrix product can be represented as the following sum of outer products:
\begin{align*}
    LU = \sum\limits_{k=1}^{n} \bm{l}^{\rm col}_k \otimes \bm{u}^{\rm row}_k.
\end{align*}
In order for $LU$ to be an integer matrix, we must have that each of $\bm{l}^{\rm col}_k \otimes \bm{u}^{\rm row}_k$ is integer. Thus,
\begin{itemize}
    \item if $\bm{l}^{\rm col}_k = \frac{1}{q_k}\bm{p}_k$, we must have that $q_k$ divides the $k^{\rm th}$ row of $U$, $\bm{u}^{\rm row}_k$.
    \item if $\bm{u}^{\rm row}_k=\frac{1}{q_k}\bm{p}_k$, $q_k$ must divide the $k^{\rm th}$ column of $L$, $\bm{l}^{\rm col}_k$.
\end{itemize}
\begin{algorithm}[H]
\caption{Full-rank integer matrix}\label{alg:LU_rank_full}
\begin{algorithmic}
\Require Construct an $n\times n$ integer matrix $A$, with ${\rm rank}(A)=n$
\State (i) Create a lower triangular matrix $L$ with arbitrary rational coefficients so that $\det(L)\neq 0$.
\State (ii) Since each column of $L$ is rational, it can be expressed as $\bm{l}^{\rm col}_k= \bm{p}_k/q_k$. Identify $q_k$ as the least-common multiple of the denominators of $\bm{l}^{\rm col}_k$ for each $k$.
\State (iii) Construct the upper triangular matrix $U$ by choosing the elements of the $k^{\rm th}$ row to be multiples of $q_k$ for each $k$.
\State (iv) Create $A$ via $A=LU$.
\end{algorithmic}
\end{algorithm}

Algorithm \ref{alg:LU_rank_full} details the general approach to constructing a full-rank matrix using matrices $L$ and $U$ that are possibly rational. The simplest case, however, is to restrict the elements of $L$ and $U$ to integers. The matrix 
\begin{align}
       A_{\rm full} = \begin{bmatrix} 
    -1 & 2 & 3 \\
    -2 & 5 & 10 \\
    -3 & 10 & 24 \end{bmatrix} =
    \begin{bmatrix}
        1 & 0 & 0 \\
        2 & 1 & 0 \\
        3 & 4 & 1
    \end{bmatrix}
    \begin{bmatrix}
        -1 & 2 & 3 \\
        0 & 1 & 4 \\
        0 & 0 & -1
    \end{bmatrix}, \label{eq:full_matrix} 
\end{align}
for example, is full rank because $U$ is full rank (by convention we consider full rank $L$).

Next, we consider how to make an integer matrix that is rank-deficient. This can be simply done by removing pivots in $U$. For example, the matrix
\begin{align}
    A_{\rm def} = \begin{bmatrix}
        1 & 2 & 3 \\
        2 & 5 & 7 \\
        3 & 10 & 13
    \end{bmatrix} = 
    \begin{bmatrix}
        1 & 0 & 0 \\
        2 & 1 & 0 \\
        3 & 4 & 1
    \end{bmatrix} 
    \begin{bmatrix}
        1 & 2 & 3 \\
        0 & 1 & 1 \\
        0 & 0 & 0
    \end{bmatrix}  \label{eq:def_matrix}
\end{align}
has $rank(A_{\rm def})=2$ since $rank(U)=2$. Moreover, the linear dependence of the columns of $A_{\rm def}$ are completely determined by that of $U$ (the third columns of $A_{\rm def}$ and $U$ are the sum of the first two of each). Algorithm \ref{alg:LU_rank_deficient} details the general approach to constructing a rank-deficient matrix.

\begin{algorithm}[H]
\caption{Constructing a rank-deficient integer matrix, $A_{\rm def}$ via $LU$-Decomposition}\label{alg:LU_rank_deficient}
\begin{algorithmic}
\Require An $n\times n$ integer matrix $A_{\rm def}$, with ${\rm rank}(A)=r<n$
\State (i) Create a lower triangular matrix $L$ with arbitrary integer coefficients.
\State (ii) Construct $r$ linearly independent columns of the upper triangular matrix $U$ by placing arbitrary integers in the off-diagonal integer coefficients.
\State (iii) Generate the remaining $n-r$ linearly dependent columns of $U$ by taking linear combinations of the others.
\State (iv) Create $A_{\rm def}$ via $A_{\rm def}=LU$.
\end{algorithmic}
\end{algorithm}

First we create a lower-triangular matrix, $L$, with arbitrary off-diagonal integer coefficients; Larger numbers would increase the arithmetical challenge. Next we choose the rank of the matrix, $r$, and generate that many linearly independent columns via Algorithm \ref{alg:LU_rank_deficient}. In step (iii) we generate the remaining $n-r$ columns as linear integer combinations of the remaining columns, guaranteeing rank-deficiency. Finally, in step (iv) we generate the original matrix multiplication, $A=LU$.

 To summarize this section, we have provided algorithms for generating both a full rank and rank-deficient matrix matrix. The $3\times 3$ matrices above ($A_{\rm full}$ from Eq.~\eqref{eq:full_matrix} and $A_{\rm def}$ from Eq.~\eqref{eq:def_matrix}) show the simplicity of the algorithm.
The main message is that the rank of the matrix can be controlled by choosing the number of non-zero diagonal elements of $U$; $rank(A_{\rm def})=2$ since $l_{33}=0$ and $rank(A_{\rm full})=3$ since $U$ is full-rank. Complications due to fractions can be circumvented with this algorithm, but the simplest case is to use integers as above.

\subsubsection{LU Decomposition and Matrix Invertibility}
 \label{sect:matrix_inversion}
Now that we know how to easily create matrices of a desired rank, we shift our focus on the ability to easily control the determinant of the matrix and subsequently the level of numerical complexity involved in the inverse of a matrix. Regarding the latter case, how do we create an invertible square matrix, $A$, so that its inverse, $A^{-1}$, only contains integers? We ask this question with the idea that students would be finding the inverse via row-reduction (which is simpler when both $A$ and $A^{-1}$ have only integers) yet we utilize the adjoint to formulate the algorithm. Suppose $A$ is a square invertible $n \times n$ matrix. Then, $\det(A)\neq 0$ and the inverse can be written as
\begin{align*}
A^{-1} = \frac{1}{\det(A)} {\rm adj}(A),
\end{align*}
where ${\rm adj}(A)$ is the adjoint of $A$. Then, for $A^{-1}$ to contain only integers, we must have that $\det(A)=\pm 1$ or $\det(A)$ divides all of the components of ${\rm adj}(A)$. Since the latter case is harder to construct, we will focus on the former. These class of matrices are known as \textbf{Unimodular matrices}. The remaining question, therefore, is to determine how to construct Unimodular matrices. The work done by \cite{Hanson1982IntegerMW} details a systematic method for constructing such a matrix, which relies on numerous row operations on the desired matrix inverse.  An easier approach would be to utilize ``special'' matrices, such as the orthogonal matrix $A$, where $A^{-1}=A^T$ and $\det(A)=\pm 1$, but significantly reduce the set of permissible matrices for our problem. We focus on utilizing the LU factorization to reach our desired result. The primary advantage of the factorization $A=LU$ is that controlling the determinant of $A$ is equivalent to controlling that of $L$ and $U$, since $\det(A) = \det(L) \det(U)$, and is relatively simple since both $L$ and $U$ are triangular. A full rank matrix such as $A_{\rm full}$ in Eq.~\eqref{eq:full_matrix} is unimodular with a ``nice'' inverse of integers because of the integer diagonals in $L$ and $U$. However, this is not a restriction, and we demonstrate one such $LU$ decomposition with non-integer rational numbers along the diagonal. Consider the following full rank matrix:

\begin{align}
    B = \begin{bmatrix}
        -1 & 2 & 3 \\
        -2 & 5 & 14 \\
        -3 & 8 & 24
    \end{bmatrix} =\begin{bmatrix}
        1 & 0 & 0 \\
        2 & 2 & 0 \\
        3 & 4 & 1
    \end{bmatrix}
    \begin{bmatrix}
        -1 & 2 & 3 \\
        0 & 1/2 & 4 \\
        0 & 0 & -1
    \end{bmatrix}.  \label{eq:lu2}
\end{align}
Both the matrices $A_{\rm full}$ and $B$ in Eqs.~\eqref{eq:full_matrix} and \eqref{eq:lu2} contain a determinant equal to 1 but their LU-decompositions are slightly different. In Eq.~\eqref{eq:full_matrix}, $\det(L)=\det(U)=1$ by only using $\pm 1$ along the diagonal. In this case, any other integers can be placed along the respective off-diagonal and the matrix $A$ will only contain integers. Slightly different is the factorization in $B$. Note that the $u_{22}$ component was relaxed to a non-integer rational number. This change in the $\det(U)=1/2$ can be offset by demanding $\det(L)=2$. Moreover, the goal is to have a matrix $B$ with only integers, so we must have that the $l_{22}$ and $l_{32}$ terms are multiples of $(1/2)^{-1}$. In general, any placement of a rational number $z$ in location $u_{ii}$ needs to be compensated by $z^{-1}$ dividing all of the elements in column $i$ of $L$. Algorithm \ref{alg:inverses} details the method of constructing Unimodular matrices and Eq.~\eqref{eq:inverse_displayed} shows the corresponding operations performed on $A_{\rm full}$ from Eq.\eqref{eq:full_matrix}.


\begin{algorithm}[H]
\caption{Constructing $A \in \mathbb{Z}^{n\times n}$ so that $A^{-1} \in \mathbb{Z}^{n\times n}$}\label{alg:inverses}
\begin{algorithmic}
\Require $A$ and $A^{-1}$ to only contain integers
\State (i) Choose the diagonal of $U$: $u_{ii} \in \mathbb{Q}$, $i=1,\ldots, n$.
\State (ii) If $u_{ii}=\pm 1$, then replace $l_{ii}$ by $\pm 1$.
\State (iii) If $u_{ii}\neq \pm 1$, replace $l_{ii}$ by $1/(u_{ii})$ and then fill the remaining column of $L$ with multiples of $l_{ii}$.
\State (iv) Fill in the remainder of $L$ and $U$ with integer coefficients.
\State (v) Generate $A$ with $\det(A)=\pm 1$ via $A=LU$. Finding the inverse via row-reduction involves no fractions.
\end{algorithmic}
\end{algorithm}

\begin{align}
   \left[ A_{\rm full} | I \right] &= \begin{bmatrix}[ccc|ccc]
    -1 & 2 & 3 & 1 & 0 & 0 \\
    -2 & 5 & 10 & 0 & 1 & 0 \\
    -3 & 10 & 24 & 0 & 1 & 0 
    \end{bmatrix} \sim
    \begin{bmatrix}[ccc|ccc]
        -1 & 2 & 3 & 1 & 0 & 0 \\
        0 & 1 & 4 & -2 & 1 & 0 \\
        0 & 0 & -1 & 5 & -4 & 1
    \end{bmatrix} \label{eq:inverse_displayed} \\ &\sim
    \begin{bmatrix}[ccc|ccc]
        -1 & 2 & 0 & 16 & -12 & 3 \\
        0 & 1 & 0 & 18 & -15 & 4 \\
        0 & 0 & -1 & 5 & -4 & 1
    \end{bmatrix} \sim
    \begin{bmatrix}[ccc|ccc]
        1 & 0 & 0 & 20 & -18 & 5 \\
        0 & 1 & 0 & 18 & -15 & 4 \\
        0 & 0 & 1 & -5 & 4 & -1
    \end{bmatrix} = \left[ I | A^{-1} \right]. \nonumber
\end{align}

In (i) we may choose the diagonal elements of $U$ to be rational. In steps (ii) and (iii) we respond to the choice in (i) by either choosing $l_{ii}=\pm 1$ or ensuring $1/u_{ii}$ divides the ith column of $L$. In (iv) we randomize the off-diagonal non-zero entries of $L$ and $U$ with integers to ensure a level of uniqueness in the problem generation. Finally, step (v) generates the unimodular matrix $A$ via $A=LU$. Algorithm \ref{alg:inverses} generates an invertible integer matrix, so that students can circumvent ugly fractions in their derivation of the inverse via row-reduction. 

To conclude this section, the $LU$ factorization can be often utilized to both control the rank of a matrix and making ``nice" integer inverses. It should be noted that, although we do not cover them here in this manuscript, the easy generation of nice inverses can be exploited in numerous subsequent factorizations. For example, instructors are often interested in providing an integer matrix that also has a nice diagonalization, $A=PDP^{-1}$, and relies on $P^{-1}$ being integer. The $LU$ factorization chooses the eigenvectors so that $P^{-1}$ is integer and combined with a choice of eigenvalues make a diagonalization and integer matrix $A$. The $LU$ factorization is used in numerous other examples. Next, in Sect.~\ref{sect:QR}, we show how it can be used to control the $QR$ factorization.

\subsubsection{QR Factorization}\label{sect:QR}
In a traditional Linear Algebra class, the students are introduced to the Gram-Schmidt process, by which a set of orthonormal vectors are obtained from linearly independent ones. Those orthonormal vectors are then used in the well-known QR factorization, $A=QR$, where $Q$ is an orthogonal matrix and $R$ is upper triangular. Since the Gram-Schmidt derived vectors are unit vectors, the matrix $Q$ often has ``ugly'' radicals. How does the clever instructor prescribe a $QR$ decomposition problem? The naive approach would be to choose an orthogonal $Q$ and random $R$, combining them together and obtaining the matrix $A$. This approach, however, would often lead to radicals in the original matrix $A$ due to the nature of the coefficients in $Q$. 
The approach we adopt here is a little different. We suggest that the instructor will want to control both the original matrix $A$ and orthogonal matrix $Q$ while using the coefficients in $R$ to connect the two together, with the singular constraint that $R$ must be upper-triangular. Therefore we ask: {\bf How do we construct an integer full-rank matrix $A$ and orthogonal matrix $Q$ independently so that $R$ is upper-triangular?} We can utilize the $LU$ factorization used in the previous section to construct a full rank matrix and the orthogonal matrix can be chosen in a number of simple ways. We proceed with a $2\times 2$ example. Suppose we choose the following:
\begin{align}
    A = LU = \begin{bmatrix}
        1 & 0 \\
        c & 1
    \end{bmatrix}
    \begin{bmatrix}
        a & b \\
        0 & d
    \end{bmatrix} =
    \begin{bmatrix}
    a & b \\
    ac & bc+d
\end{bmatrix},
\end{align}
where $a, b, c$ and $d$ are free parameters that we can choose, representing 4 degrees of freedom with the constraint that $a,d\neq 0$ from the linear independence of the columns of $A$. Then, if $A=QR$, by the orthogonality of $Q$ we have
\begin{align*}
    R =  Q^T A = \begin{bmatrix}
        q_{11} & q_{21} \\
        q_{12} & q_{22}
    \end{bmatrix}
    \begin{bmatrix}
        a & b \\
        ac & bc+d
    \end{bmatrix} = \begin{bmatrix}
        a q_{11} + ac q_{21} & b q_{11} + q_{21}(bc+d) \\
        a q_{12} + ac q_{22} & bq_{12} + q_{22} (bc+d)
    \end{bmatrix}
\end{align*}
Now the only constraint that we have on $R$ is that it is upper-triangular, which imposes the condition on the 3rd component:
\begin{align}
    r_{12}=a\left( q_{12} + c q_{22} \right)= 0. \label{eq:R_condition}
\end{align}
Equation \eqref{eq:R_condition} simply states that the second column of $Q$ must be orthogonal to the first column of $A$, which is how Gram-Schmidt works. In order to reduce the size of the parameter space we assume that $a, b, d, q_{11}, q_{12}, q_{21},$ and $q_{22}$ are initially chosen and that $c$ is determined via Eq.~\eqref{eq:R_condition} as $c=-q_{12}/q_{22}$ when $q_{22}\neq 0$ (since $a\neq 0$). This includes the case $q_{12}=0$, where we have $Q=I$ so that $c=0$ and any full rank upper triangular matrix $A$ would suffice, resulting in $R=Q^T A = IA = A$. The case $q_{22}=0$ requires the $A=PLU$ factorization, since $a\neq 0$, so we omit it for simplicity. Therefore, for the algorithm we consider two cases:
\begin{enumerate}[(a)]
    \item If $Q=I$, then, choose any matrix of the form $A=\begin{bmatrix}
     a & b \\ 0 & d
    \end{bmatrix}$, where $a,b,d \in \mathbb{Z}$ and $a,d \neq 0$. Then, generate $R$ via $R=Q^T A = A$.
    \item If $Q\neq I$, then set $c=-q_{12}/q_{22}$ and generate $A$ via $
        A = LU =\begin{bmatrix}
            1 & 0 \\ c & 1
        \end{bmatrix}
        \begin{bmatrix}
            a & b \\
            0 & d
        \end{bmatrix}$. Then, generate $R$ via $R=Q^T A$. 
\end{enumerate}

\begin{algorithm}[H]
\caption{Constructing $A=QR$ where $A$ is an integer matrix and $Q$ is orthogonal}\label{alg:QR}
\begin{algorithmic}
\Require Integer matrix, $A$, Orthogonal matrix $Q$ and Upper-triangular $R$.
\State (i) Choose numbers for an upper triangular matrix $U$ with $u_{ii}\neq 0$ for $i=1,2,\ldots, n$.
\State (ii) Parameterize the non-zero off-diagonal components in $L$ and choose the diagonal elements to be 1.
\State (iii) Choose a vector $\bm{v}\in\mathbb{R}^n$ that has integer components. Generate a non-identity orthogonal matrix $Q$ via the Householder transformation $Q=I-2\bm{v}\bm{v}^T/||\bm{v}||^2$. The resultant matrix $Q$ will then be rational (possibly integer).
\State (iv) Find the values of the $L$ parameterized variables that would make the entries below the diagonal of $R=Q^T A$ equal to 0.
\State (v) Find $A$ via $A=c LU$ and $R$ via $R=Q^T A$, where $c$ (the lcd of the components of $A$) is chosen to make $A$ integer.
\end{algorithmic}
\end{algorithm}

\begin{align}
&{}    &&U= \begin{bmatrix}
        1 & 2 & 3 \\
        0 & 4 & 5 \\
        0 & 0 & 6
    \end{bmatrix}, \quad
    L = \begin{bmatrix}
        1 & 0 & 0 \\
        x & 1 & 0 \\
        y & z & 1
    \end{bmatrix}, \label{eq:QR_3by3_example1} \\
&(\bm{v} \in \mathbb{Z}^3, ||\bm{v}||=1), \; &&\bm{v}= \begin{bmatrix}
        1 \\ 0 \\ 0
    \end{bmatrix}, \quad Q = \begin{bmatrix}
        -1 & 0 & 0 \\
        0 & 1 & 0  \\
        0 & 0 & 1
    \end{bmatrix}, \quad A = \begin{bmatrix}
        1 & 2 & 3 \\
        0 & 4 & 5 \\
        0 & 0 & 6
    \end{bmatrix}, \label{eq:QR_3by3_example2} \\
   &(\bm{v} \in \mathbb{Z}^3, ||\bm{v}||\neq 1,2), \; &&\bm{v}= \begin{bmatrix}
        1 \\ 2 \\ 3
    \end{bmatrix}, \quad Q = \frac{1}{7}\begin{bmatrix}
    6 & -2 & -3 \\
    -2 & 3 & -6 \\
    -3 & -6 & -2
    \end{bmatrix}, \quad A = 6\begin{bmatrix}
    1 & 2 & 3 \\
    \tfrac{-1}{3} & \tfrac{10}{3} & 4 \\
    \tfrac{-1}{2} & -13 & \tfrac{-21}{2}
    \end{bmatrix}. \label{eq:QR_3by3_example4}
\end{align}

Algorithm~\ref{alg:QR} details the method of constructing a $QR$ decomposition that yields an integer matrix and Eqs.\eqref{eq:QR_3by3_example1}--\eqref{eq:QR_3by3_example4} provide a few $3 \times 3$ examples with key characteristics. In (i) we first choose a full-rank upper-triangular matrix $U$ with random integer coefficients. In (ii) we parameterize the non-zero off-diagonal components of the upper-triangular matrix $L$ and impose that $L$ is {\bf unit triangular} ($1$'s along the diagonal). For each of the examples above we utilize $U$ and $L$ from Eq.~\eqref{eq:QR_3by3_example1}. 

In (iii) we generate an $n\times n$ orthogonal matrix $Q$ via the Householder transformation $Q=I-2\bm{v}\bm{v}^T/||\bm{v}||^2$, where the vector $\bm{v}\in\mathbb{R}^n$ is chosen with integer components. 
The ``niceness'' of $Q$, then, relies on the choice of $\bm{v}$. If, for example, we choose integer $\bm{v}$ with $||\bm{v}||=1$ (as in Eq.~\eqref{eq:QR_3by3_example2}), the matrix $Q$ will be entirely integer (since $\bm{v}\bm{v}^T$ will be integer). For an n-dimensional vector this amounts to $2n$ choices (e.g. in 3D  $\begin{bmatrix}
    \pm 1 & 0 & 0
\end{bmatrix}^T, \begin{bmatrix}
    0 & \pm 1 & 0
\end{bmatrix}^T, \begin{bmatrix}
    0 & 0 & \pm 1
\end{bmatrix}^T$). Another way to generate an integer orthogonal matrix $Q$ is to use $||\bm{v}||^2=2$, but Algorithm \ref{alg:QR} would need to incorporate $PLU$ factorization and is not considered here.

For randomly chosen integer $\bm{v}$, such as that in Eq.~\eqref{eq:QR_3by3_example4} the matrix $Q$ will be entirely rational. It should be noted that any rational $\bm{v}_1$ can be made into a integer vector $\bm{v}_2$ that gives the same orthogonal matrix $Q$. For this reason, and for simplicity, we consider only integer $\bm{v}$. 

In (iv) the parameterized values of $L$ are then found by solving orthogonality relationships (each column of $Q$ must be orthogonal to the preceding columns of $A$). Finally, in (v) $A$ is formed by scaling $LU$ so that $A=cLU$ is integer and $R$ is generated via $R=Q^T A$. 

The ingredient in this algorithm that binds $A$ to $Q$ is the constraint that the lower-diagonal of $R$ must be zero. In some cases, for example in Eq.~\eqref{eq:QR_3by3_example2}), this simply leads to $x=y=z=0$. In others, the solution is more arithmetically complex (e.g. in Eq.~\eqref{eq:QR_3by3_example4}, $x=-1/3, y=-1/2$ and $z$=-3).  The orthogonality relationships for a $3\times 3$ matrix $R$ are equivalent to
\begin{align*}
    r_{21} = c\bm{q}_2^T \cdot \begin{bmatrix}
        \bm{l}^{\rm row}_1 \cdot \bm{u}^{\rm col}_1 \\
        \bm{l}^{\rm row}_2 \cdot \bm{u}^{\rm col}_1 \\
        \bm{l}^{\rm row}_3 \cdot \bm{u}^{\rm col}_1
    \end{bmatrix} = 0, \quad
    r_{31} = c\bm{q}_3^T \cdot \begin{bmatrix}
        \bm{l}^{\rm row}_1 \cdot \bm{u}^{\rm col}_1 \\
        \bm{l}^{\rm row}_2 \cdot \bm{u}^{\rm col}_1 \\
        \bm{l}^{\rm row}_3 \cdot \bm{u}^{\rm col}_1
    \end{bmatrix} = 0, \quad
    r_{32} = c\bm{q}_3^T \cdot \begin{bmatrix}
    \bm{l}^{\rm row}_1 \cdot \bm{u}^{\rm col}_2 \\
    \bm{l}^{\rm row}_2 \cdot \bm{u}^{\rm col}_2 \\
    \bm{l}^{\rm row}_3 \cdot \bm{u}^{\rm col}_2
    \end{bmatrix} = 0,
\end{align*}
where $\bm{q}_1, \bm{q}_2, \bm{q}_3$ represent the columns of $Q$, $\bm{l}^{\rm row}_1, \bm{l}^{\rm row}_2, \bm{l}^{\rm row}_3$ represent the rows of $L$, and $\bm{u}^{\rm col}_1, \bm{u}^{\rm col}_2, \bm{u}^{\rm col}_3$ represent the columns of $U$ and $R = c Q^T L U= c \begin{bmatrix}
        \bm{q}_1 & \bm{q}_2 & \bm{q}_3
    \end{bmatrix}^T \begin{bmatrix}
        \bm{l}^{\rm row}_1 &
        \bm{l}^{\rm row}_2 &
        \bm{l}^{\rm row}_3
    \end{bmatrix}^T \begin{bmatrix}
        \bm{u}^{\rm col}_1 & \bm{u}^{\rm col}_2 & \bm{u}^{\rm col}_3
    \end{bmatrix}.$  
In general, the orthogonality relationships are
\begin{align}
    r_{ij} = c\bm{q}_i^T \cdot \begin{bmatrix}
        \bm{l}^{\rm row}_1 \cdot \bm{u}^{\rm col}_j \\
        \bm{l}^{\rm row}_2 \cdot \bm{u}^{\rm col}_j \\
        \vdots \\
        \bm{l}^{\rm row}_n \cdot \bm{u}^{\rm col}_j
    \end{bmatrix} = 0, \qquad j=1, 2, \ldots, n, \quad i>j, \label{eq:R_relationship_general}
\end{align}
and can be solved easily using a computer algebra system (such as Mathematica). It should be noted that the consistency of the system given by Eq.~\eqref{eq:R_relationship_general} is not guaranteed. For example, the case of $\bm{v}\in \mathbb{Z}^n$  and $||v||^2=2$ would yield such an inconsistency and require that Algorithm \ref{alg:QR} be modified with the $PLU$ factorization of $A$. Since, the user has control over the components of $\bm{v}$ we consider such a restriction unnecessary as the components can be simply changed. There are also ways of making integer $Q$ with rational numbers. For example, the vector $\bm{v} = \begin{bmatrix}
    \sqrt{2} & 0 & 0
\end{bmatrix}^T$ has $||\bm{v}||=2$ and generates an integer $Q$. The irrational complexity can be obviously increased to yield arithmetically more difficult $Q$, but we omit these cases from Algorithm \ref{alg:QR}, for simplicity.

It should be noted that other algorithms can be used to generate an orthogonal matrix. In 2D, the matrix $Q=\frac{1}{a^2+b^2}\begin{bmatrix}
    a & -b \\
    b & a
\end{bmatrix}$ is orthogonal and can be easily generated. In 3D, we can start with orthogonal vectors (which is non-trivial) and generate a third via the cross-product. The columns of $Q$ are then the unit versions of those choices. In higher dimensions we could utilize the Gram-Schmidt algorithm, but this defeats the purpose of controlling the ``niceness'' of $Q$. We utilize the Householder transformation idea above since all of other methods are not generalizable or limit the user control. Finally, note that this algorithm will give you an orthogonal matrix $Q$ that is consistent with Gram-Schmidt up to a constant multiple of a column (since orthogonality is not unique).

\subsection{Constructing a system of equations} \label{sect:linear_systems}
Now that we know how to control the properties of a matrix, we shift our focus to setting up a system of equations. How does one create a system of equations, represented by $A\bm{x}=\bm{b}$, with a unique solution? To answer this question we note the following: (1) a solution exists only if $\bm{b}$ can be written as a linear combination of the columns of $A$ and (2) the columns must be linearly independent for uniqueness. Algorithm \ref{alg:lin_alg1} outlines the process. We proceed by example.

\begin{algorithm}
\caption{Constructing $A\bm{x}=\bm{b}$ where $\bm{x}$ is unique and $\bm{b}$ is unspecified}\label{alg:lin_alg1}
\begin{algorithmic}
\Require Generate an integer matrix $A$, vector $\bm{b}$ and unique solution $\bm{x}$.
\State (i) Generate linearly independent column vectors of $A$, $\bm{v}_1, \bm{v}_2, \ldots \bm{v}_n$ where $\bm{v}_j \in \mathbb{R}^n$, from $A=LU$.
\State (ii) Choose an integer solution $\bm{x}=[x_1, x_2, \ldots, x_n]^T$.
\State (iii) Create $\bm{b}$ as a linear combination of the column vectors $\bm{b}=x_1 \bm{v}_1 + x_2 \bm{v}_2 + \ldots x_n \bm{v}_n$
\end{algorithmic}
\end{algorithm}

In step (i), we create a full-rank matrix $A$, via Algorithm \ref{alg:LU_rank_full}.
\begin{align}
    A &= \begin{bmatrix}
        1 & 0 & 1 \\
        -1 & 1 & 2 \\
        2 & 2 & 3
    \end{bmatrix}=\begin{bmatrix}
        1 & 0 & 0 \\
        -1 & 1 & 0 \\
        2 & 2 & 1
    \end{bmatrix} 
    \begin{bmatrix}
        1 & 0 & 1 \\
        0 & 1 & 3 \\
        0 & 0 & -5
    \end{bmatrix} = [\bm{v}_1, \bm{v}_2, \bm{v}_3],  \label{eq:lin_system1_example1} \\
    \bm{x} &= [1, 2, 3]^T,   \label{eq:lin_system1_example2} \\
    \bm{b} &= 1\bm{v}_1 + 2\bm{v}_2 + 3\bm{v}_3.   \label{eq:lin_system1_example3}
\end{align}

Since both $L$ and $U$ in Eq.~\eqref{eq:lin_system1_example1} have a full set of pivots, so does $A$, ensuring linear independence of the columns, $\bm{v}_1, \bm{v}_2$, and $\bm{v}_3$. In step (ii) (Eq.~\eqref{eq:lin_system1_example2}) an integer solution $x_1, x_2$, and $x_3$ is chosen. Finally, in (iii) the vector $\bm{b}$ is constructed as a linear combination of the columns (Eq.~\eqref{eq:lin_system1_example3}), ensuring existence. Of note, is that we needed the elements of $L$, $U$, and $\bm{x}$ to be integers to ensure that $A$ does as well. Similar to Algorithm \ref{alg:inverses} above, we can relax this restriction to rational numbers in, by choosing the numbers in $L$ to account for the rational numbers in $U$ and vice versa. 

Algorithm \ref{alg:lin_alg1} gives the user control over the coefficients of the matrix and the solution while ensuring $\bm{b}\in C(A)$ (the column space of $A$). Although this construction methodology seems trivial given the ``column picture'' of a system of linear equations, the difficulty lies in the order of the algorithm. Suppose, for example, that we have chosen the column vectors $\bm{v}_i$ and the vector $\bm{b}$. In that case we could not guarantee that the solution $\bm{x}$ only contains integers. Suppose instead that we wanted to control $\bm{x}$ and $\bm{b}$. We couldn't guarantee, then, that the coefficients of the matrix are integers. For example, suppose that $A$ is $2\times 2$. Choosing $\bm{x}$ and $\bm{b}$ would yield a 2 degree-of-system (4 unknowns and 2 equations) with additional restrictions on the coefficients needed. It is crucial therefore, to follow the flow of Algorithm \ref{alg:lin_alg1}. 

Suppose instead that we wanted to construct a system with infinitely many solutions. How could we construct such a system? The easiest way is to follow Algorithm \ref{alg:lin_alg1} but choose at least one of the columns of $A$ to depend on the others.

\begin{algorithm}
\caption{Constructing $A\bm{x}=\bm{b}$ where $\bm{x}$ is not unique} \label{alg:infinite_soln}
\begin{algorithmic}
\Require Integer coefficients in $A$ and $\bm{b}$.
\State (i) Choose $n-1$ integers to represent a linear combination of the columns.
\State (ii) Randomly choose integer elements for the lower-diagonal of $L$ and non-zero components of the first $n-1$ columns of $U$.
\State (iii) Determine the final column $\bm{v}_n$ via the linear combination chosen in (i).
\State (iv) Generate the matrix $A$ containing $n-1$ linearly independent vectors $\bm{v}_1, \bm{v}_2, \ldots \bm{v}_{n-1}$ and one linearly dependent vector $\bm{v}_n$ via $A=LU$. The chosen column dependence of $U$ is the same as $A$.
\State (v) Choose a particular solution $x_1, x_2, \ldots, x_n$ so that $\bm{b} = x_1 \bm{v}_1 + x_2 \bm{v}_2 + \cdots x_n \bm{v}_n$.
\end{algorithmic}
\end{algorithm}

We proceed with an $3 \times 3$ example for Algorithm \ref{alg:infinite_soln}. In step (i), we choose an integer column dependence (i.e. $\bm{u}_3 = c_1 \bm{u}_1 + c_2 \bm{u}_2$). In Eq.\eqref{eq:col_dependence} we have chosen integers $\underline{1}$ and $\underline{2}$. In step (ii)-(iv) (Eq.~\eqref{eq:lin_system2_example1}), we use an adaptation of Algorithm \ref{alg:LU_rank_deficient} to generate a rank-deficient matrix $A_{\rm def}=LU$ whose columns, $\bm{v}_1, \bm{v}_2, \bm{v}_3$, form a linearly dependent set. In particular, the lower-diagonal of $L$ (colored {\color{red} red}), and the non-zero components of the first two columns of $U$ ({\color{red} red}) are randomly chosen amongst $\mathbb{Z}$ in step (ii). In step (iii) the linear combination chosen in (i) is used to generate the final column $\bm{u}_3$ of $U$ and subsequently the rank-deficient matrix is formed in (iv) via $A_{\rm def}=LU$. The column vectors given in Eq.~\eqref{eq:lin_system1_example1}, $\bm{v}_1, \bm{v}_2, \bm{v}_3$ have the same linear combination as those in $U$, namely $\bm{v}_3 = \underline{1} \bm{v}_1 + \underline{2}\bm{v}_2$. This removes the uniqueness of a solution to $A\bm{x}=\bm{b}$. In step (v) an integer solution $x_1$, $x_2$, and $x_3$ is chosen as before, and a vector $\bm{b}$ is constructed as a linear combination of the columns (Eq.~\eqref{eq:lin_system2_example5}), ensuring existence. Note that in Eq.~\eqref{eq:lin_system2_example6} the row reduction of the augmented form leads to the desired result. Of particular note, the chosen column dependence shows up in the reduced-row echelon form of $A$.

\begin{align}
    \bm{u}_3 &= {\color{blue} \underline{1}} \bm{u}_2 + {\color{blue} \underline{2}} \bm{u}_3 \label{eq:col_dependence}\\
    A_{\rm def} &= LU=\begin{bmatrix}
        1 & 0 & 0 \\
        {\color{red} -1} & 1 & 0 \\
        {\color{red} 2} & {\color{red} 2} & 1
    \end{bmatrix}     
\begin{bmatrix}
        {\color{red} 2} & {\color{red} 4} & 10\\
        0 & {\color{red} 3} & 6 \\
        0 & 0 & 0 
    \end{bmatrix}= \begin{bmatrix}
        2 & 4 & 10 \\
        -2 & -1 & -4 \\
        4 & 14 & 32
    \end{bmatrix} = \begin{bmatrix}
        \vert & \vert & \vert \\
        \bm{v}_1 & \bm{v}_2 & \bm{v}_3 \\
        \vert & \vert & \vert
    \end{bmatrix}
\label{eq:lin_system2_example1}\\
    \bm{x} &= [1, 2, 3]^T,   \label{eq:lin_system2_example4} \\
    \bm{b} &= 1\bm{v}_1 + 2\bm{v}_2 + 3\bm{v}_3 = \begin{bmatrix} 40 & -16 & 128 \end{bmatrix}^T  \label{eq:lin_system2_example5}\\
    [A_{\rm def} |\, \bm{b}\,] &= \begin{amatrix}{3}
        2 & 4 & 10 & 40 \\
        -2 & -1 & -4 & -16 \\
        4 & 14 & 32 & 128
    \end{amatrix} = \begin{amatrix}{3}
        1 & 0 & {\color{blue} \underline{1}} & 4 \\
        0 & 1 & {\color{blue} \underline{2}} & 8 \\
        0 & 0 & 0 & 0
    \end{amatrix}
    \label{eq:lin_system2_example6}
\end{align}

The main idea behind Algorithm \ref{alg:infinite_soln} is to use linear dependence of the columns to remove uniqueness while maintaining existence by choosing $b \in C(A)$. In the above example, we chose the last column to depend on the first. This process can be generalized by taking $n-k$ linearly independent ones and choosing $k$ linear combinations of the remaining dependent ones as in Eq.~\eqref{eq:col_dependence}. It should be noted that choosing the dependence of the columns of $U$ is imperative for the algorithm to work. Even random integer coefficients in $U$ may lead to integer $A$ that row reduces to fractional reduced-row echelon form, which is less desirable from an educational point of view.

Finally, how do we create a system of linear equations that has no solution? This question, is actually more difficult than the previous questions. We will rely on the orthogonality of the four fundamental spaces. Algorithm \ref{alg:no_soln} details the process.

\begin{algorithm}
\caption{Constructing $A\bm{x}=\bm{b}$ where $\bm{x}$ does not exist} \label{alg:no_soln}
\begin{algorithmic}
\Require No solution to $A\bm{x} = \bm{b}$ where both $A$ and $\bm{b}$ contain integers.
\State (i) Choose $n-1$ integers to represent a linear combination of the columns.
\State (ii) Follow the procedure in Algorithm \ref{alg:LU_rank_deficient} to generate a rank-deficient matrix, $A_{\rm def}=LU$, but make the last column of $U$ the above integer combination of the preceding columns.
\State (iii) Generate the matrix $A$ containing $n-1$ linearly independent vectors $\bm{v}_1, \bm{v}_2, \ldots \bm{v}_{n-1}$ and one linearly dependent vector $\bm{v}_n$ via $A=LU$. The chosen column dependence of $U$ is the same as $A$.
\State (iv) Find an element in the left nullspace, $\bm{y} \in N(A^T)$: $A^T \bm{y} = \bm{0}$. 
\State (v) Choose random integer coefficients $x_1, x_2, \ldots, x_n$ and let $\bm{b} = x_1 \bm{v}_1 + x_2 \bm{v}_2 + \cdots x_n \bm{v}_n + \bm{y}$ so that $\bm{b}$ is orthogonal to the span of the columns $\{\bm{v}_1, \bm{v}_2, \ldots, \bm{v}_{n} \}$. 
\State (vi) Create $A$ using the column vectors $\bm{v}_1, \bm{v}_2, \ldots, \bm{v}_n$.
\end{algorithmic}
\end{algorithm}


We proceed with an $3 \times 3$ example for Algorithm \ref{alg:no_soln} and use the same matrix, $A_{\rm def}$, and vectors, $\bm{v}_1, \bm{v}_2, \bm{v}_3$, used for the example in Algorithm \ref{alg:infinite_soln}. At the conclusion of step (iii), we have two linearly independent vectors $\bm{v}_1$ and $\bm{v}_2$ along with a third vector $\bm{v}_3$, which is linearly dependent ($\bm{v}_1+\underline{1}\bm{v}_2=\underline{2}\bm{v}_3$). In step (iv), we find an element in the left nullspace of $A_{\rm def}$. For example, it can be verified in Eq.~\eqref{eq:lin_system3_example4} that $\begin{bmatrix} -4 & -2 & 1 \end{bmatrix}^T$ is in $N(A^T)$. In (iv) an integer solution $x_1$, $x_2$, and $x_3$ is chosen as before, and a vector $\bm{b}$ is constructed as a linear combination of the columns plus the y vector (Eq.~\eqref{eq:lin_system3_example6}) so that $\bm{b}$ is not in $C(A)$. Note that the row reduction of $\begin{bmatrix} A|\bm{b} \end{bmatrix}$ produces the desired inconsistency as shown in Eq.~\eqref{eq:lin_system3_example7}.

\begin{align}
    A_{\rm def}^T \bm{y} &= \begin{bmatrix}
        2 & 4 & 10 \\
        -2 & -1 & -4 \\
        4 & 14 & 32
    \end{bmatrix} \begin{bmatrix} -4 \\ -2 \\ 1 \end{bmatrix} = \begin{bmatrix} 0 \\ 0 \\ 0 \end{bmatrix}\label{eq:lin_system3_example4}
    \\
    \bm{x} &= [1, 2, 3]^T,   \label{eq:lin_system3_example5} \\
    \bm{b} &= 1\bm{v}_1 + 2\bm{v}_2 + 3\bm{v}_3 + \bm{y}= \begin{bmatrix} 12 & -14 & 115\end{bmatrix}^T  \label{eq:lin_system3_example6}\\
    [A_{\rm def} |\, \bm{b}\,] &= \begin{amatrix}{3}
        2 & 4 & 10 & 12 \\
        -2 & -1 & -4 & -14 \\
        4 & 14 & 32 & 115
    \end{amatrix} = \begin{amatrix}{3}
        1 & 0 & {\color{blue} \underline{1}} & 0 \\
        0 & 1 & {\color{blue} \underline{2}} & 0 \\
        0 & 0 & 0 & 1
    \end{amatrix}
    \label{eq:lin_system3_example7}
\end{align}
 
Both the linearly independent set of vectors $V=\{\bm{v}_1, \bm{v}_2, \ldots, \bm{v}_{n-1}\}$ containing integer coefficients and the linear dependence $\bm{v}_n = k_1 \bm{v}_1 + \cdots k_{n-1} \bm{v}_{n-1}$ can be obtained following steps (i) and (ii) in Algorithm \ref{alg:infinite_soln}. Then since $rank(A)<n$, we can find a non-trivial vector $\bm{y}$ in the left nullspace. By orthogonality of $N(A^T)$, the left nullspace, and $C(A)$, this vector $\bm{y}$ is orthogonal to the linearly independent vectors in $V$. We then choose a linear combination and perturb $\bm{b}=x_1 \bm{v}_1 + \cdots \bm{v}_n + \bm{y}$ so that $\bm{b}$ does not lie in $C(A)$. The linear combination can be chosen as the user likes and does not affect the final outcome since $\bm{y} \in N(A^T)$. The resultant system $A\bm{x}=\bm{b}$ would then be inconsistent as $\bm{b}$ has a component not spanned by the columns of $A$.

\section{Conclusion}

As an instructor, designing lessons that are both accessible and have a reasonably ``nice'' (integer) solution is far from trivial. Typically, educators rely on course textbooks, online instructional tools, or other problem-generating software, which often lead to complex answers that may hamper students' understanding. In this paper, we have proposed numerous algorithms for designing problems in the context of linear algebra with nice integer solutions. Our algorithmic approach provides educators greater control over the parameters of the problems, and can be used to enhance instruction and improve student confidence. Regarding the properties of the matrix, we proposed seven algorithms to control the rank of the matrix, the invertibility of the matrix, and the $QR$ decomposition. We used these algorithms to formulate systems of equations with unique, infinite, or no solutions. We show that by manipulating elements in a matrix, problems can be generated that have nice interim and final solutions, reducing algebraic stumbling blocks for students and keeping the primary focus on the concept being taught. Additionally, we find that the scholarly pursuit of reverse engineering for algorithm development, leveraging matrix properties and linear algebra theory, is a substantial academic undertaking that can stand on its own for both students and faculty members. 

The complete source code and algorithms presented in this article can be accessed and downloaded from the dedicated repository \url{https://github.com/Ryallaire/REVENG} to facilitate further research, collaboration, and reproducibility.




\bibliographystyle{plain}
\bibliography{mybib}

\begin{thebibliography}{10}

\bibitem{doi:10.1111/1467-8721.00196}
Mark~H. Ashcraft.
\newblock Math anxiety: Personal, educational, and cognitive consequences.
\newblock {\em Current Directions in Psychological Science}, 11(5):181--185,
  2002.

\bibitem{Baloglu}
Mustafa Baloğlu and Recep Koçak.
\newblock Multivariate investigation of the differences in mathematics anxiety.
\newblock {\em Personality and Individual Differences}, 40:1325--1335, 05 2006.

\bibitem{Bistritz.2010}
Y.~Bistritz and Yaron Segalov.
\newblock Fraction-free inversion of a toeplitz matrix.
\newblock 2010.

\bibitem{10.2307/2691448}
Jim Buddenhagen, Charles Ford, and Mike May.
\newblock Nice cubic polynomials, pythagorean triples, and the law of cosines.
\newblock {\em Mathematics Magazine}, 65(4):244--249, 1992.

\bibitem{doi:10.1080/00029890.1986.11971837}
Samuel Councilman.
\newblock Eigenvalues and eigenvectors of “n-matrices”.
\newblock {\em The American Mathematical Monthly}, 93(5):392--395, 1986.

\bibitem{doi:10.1080/07468342.2000.11974176}
Al~Cuoco.
\newblock Meta-problems in mathematics.
\newblock {\em The College Mathematics Journal}, 31(5):373--378, 2000.

\bibitem{doi:10.1080/00029890.1984.11971414}
W.~P. Galvin.
\newblock Matrices with “custom-built” eigenspaces.
\newblock {\em The American Mathematical Monthly}, 91(5):308--309, 1984.

\bibitem{10.2307/2322394}
Richard~C. Gilbert.
\newblock Companion matrices with integer entries and integer eigenvalues and
  eigenvectors.
\newblock {\em The American Mathematical Monthly}, 95(10):947--950, 1988.

\bibitem{10.2307/3026530}
Robert Hanson.
\newblock Integer matrices whose inverse contain only integers.
\newblock {\em The Two-Year College Mathematics Journal}, 13(1):18--21, 1982.

\bibitem{Hanson1982IntegerMW}
Robert Hanson.
\newblock Integer matrices whose inverse contain only integers.
\newblock {\em Two-Year College Mathematics Journal}, 13:18, 1982.

\bibitem{Hoffman2010}
Bobby Hoffman.
\newblock "i think i can, but i'm afraid to try": The role of self-efficacy
  beliefs and mathematics anxiety in mathematics problem-solving efficiency.
\newblock {\em Learning and Individual Differences}, 20(3):276--283, June 2010.

\bibitem{oro24261}
Sue Johnston-Wilder and Clare Lee.
\newblock Developing mathematical resilience.
\newblock In {\em BERA Annual Conference 2010}, 2010.

\bibitem{LEE2017269}
Clare Lee and Sue Johnston-Wilder.
\newblock Chapter 10 - the construct of mathematical resilience.
\newblock In Ulises {Xolocotzin Eligio}, editor, {\em Understanding Emotions in
  Mathematical Thinking and Learning}, pages 269--291. Academic Press, San
  Diego, 2017.

\bibitem{MARTIN200853}
Andrew~J. Martin and Herbert~W. Marsh.
\newblock Academic buoyancy: Towards an understanding of students' everyday
  academic resilience.
\newblock {\em Journal of School Psychology}, 46(1):53--83, 2008.

\bibitem{https://doi.org/10.48550/arxiv.0712.3060}
Greg Martin and Erick~B. Wong.
\newblock Almost all integer matrices have no integer eigenvalues, 2007.

\bibitem{Martin.2008}
Greg Martin and Erick~B. Wong.
\newblock The number of 2×2 integer matrices having a prescribed integer
  eigenvalue.
\newblock {\em Algebra \& Number Theory}, 2008.

\bibitem{strang2023}
G~Strang.
\newblock {\em Introduction to Linear Algebra (6th ed.)}.
\newblock Cambridge University Press, Cambridge, 2023.

\bibitem{Towse}
Christopher Towse and Eric Campbell.
\newblock Constructing integer matrices with integer eigenvalues.
\newblock {\em Applied Probability Trust}, 03 2016.

\end{thebibliography}

\end{document}